\documentclass[12pt]{article}
\usepackage{theorem,amsfonts}

\newtheorem{theorem}{Theorem}
\newtheorem{proposition}{Proposition}
\newtheorem{lemma}{Lemma}
\newtheorem{corollary}{Corollary}
\theorembodyfont{\upshape}
\newtheorem{definition}{Definition}
\newtheorem{remark}{Remark}
\newtheorem{proof}{Proof}

\newtheorem{acknowledgement}{Acknowledgement}

\newcommand{\bt}{\begin{theorem}}
\newcommand{\et}{\end{theorem}}
\newcommand{\bl}{\begin{lemma}}
\newcommand{\el}{\end{lemma}}
\newcommand{\bp}{\begin{proposition}}
\newcommand{\ep}{\end{proposition}}
\newcommand{\bo}{\begin{proof}}
\newcommand{\eo}{\end{proof}}
\newcommand{\br}{\begin{remark}}
\newcommand{\er}{\end{remark}}
\newcommand{\bd}{\begin{definition}}
\newcommand{\ed}{\end{definition}}
\newcommand{\bc}{\begin{corollary}}
\newcommand{\ec}{\end{corollary}}
\newcommand{\be}{\begin{enumerate}}
\newcommand{\ee}{\end{enumerate}}

\title{Operator semi-selfdecomposable measures and related nested
subclasses of measures on vector spaces}
\author{C. R. E. Raja}
\date{}

\begin{document}
\maketitle

\let\epsi=\epsilon
\let\vepsi=\varepsilon
\let\lam=\lambda
\let\Lam=\Lambda 
\let\ap=\alpha
\let\vp=\varphi
\let\ra=\rightarrow
\let\Ra=\Rightarrow 
\let\LRa=\Leftrightarrow
\let\Llra=\Longleftrightarrow
\let\Lla=\Longleftarrow
\let\lra=\longrightarrow
\let\Lra=\Longrightarrow
\let\ba=\beta
\let\ga=\gamma
\let\Ga=\Gamma
\let\un=\upsilon

\begin{abstract}
$T$-semi-selfdecomposability, operator selfdecomposability and subclasses 
$L_m(b, (T_t))$ and $\tilde L_m(b, (T_t))$ of measures on complete 
separable metric vector spaces are introduced and basic properties are 
proved.  In particular, we show that $\mu$ is $T$-semi-selfdecomposable if 
and only if $\mu = T(\mu ) \nu$ where $\nu$ is infinitely divisible and 
$\mu$ is operator selfdecomposable if and only if $\mu \in L_0(b, (T_t))$ 
for all $0< b < 1$.  
\end{abstract}



\begin{section}{Introduction}

Let $E$ be a complete separable metric topological vector space.  
Let ${\cal P}(E)$ be the space of all regular Borel
probability measures on $E$.  We endow ${\cal P}(E)$ with the weak
topology, a net $(\mu _i)$ in ${\cal P}(E)$ converges weakly to $\mu$
if $\mu _i (f) \ra \mu (f)$ for all bounded continuous functions 
on $E$.  Then ${\cal P}(E)$ with convolution is a complete 
separable metric topological semigroup.  For $\mu , \lam \in {\cal P}(E)$,
$\mu \lam$ denotes the convolution product of $\mu$ and $\lam$.  

A triangular system $\{ \lam _{n,j}\mid n \geq 1,~ 1\leq j \leq k_n \}$ of
probability measures on $E$ is said to be {\it infinitesimal} if 
$$ \lim _{n \ra \infty} \sup _j \lam _{n,j}
(E\setminus U) = 0 $$ for each open neighborhood $U$ of $0$.  We say that
a triangular system $(\lam _{n,j})$ converges to $\lam$ if 
$\prod _{j=1}^{k_n}\lam _{n,j} \ra \lam$.   

A probability measure $\mu$ on $E$ is
called {\it infinitely divisible} if for each $n \geq 1$, there exists a
$\mu _n \in {\cal P}(E)$ such that $$\mu _n ^n= \mu .$$ Let $I(E)$ be the
set of all infinitely divisible measures in ${\cal P}(E)$.

A continuous convolution semigroup (c.c.s.), $(\mu _t)$ of probability 
measures on $E$ is a continuous semigroup homomorphism $t\mapsto \mu _t$ 
from $[0, \infty )$ into ${\cal P}(E)$.

A probability measure $\mu$ on $E$ is called {\it embeddable} if there 
exists a c.c.s. $(\mu _t)$ such that $\mu _1 = \mu$.  It is easy to see 
that embeddable 
measures are infinitely divisible.  We will see that the converse is true
if $E$ is strongly root compact.  We recall that $E$ is said to be  
{\it strongly root compact} if for every compact set $C\subset E$ there
exists a compact set $C_0 \subset E$ such that for every $n \in \mathbb N$
all finite sequences $\{ x_1, x_2, \cdots , x_n\}$, $x_n =0$, satisfying 
the condition $$(C+x_i +C+x_j) \cap (C+x_{i+j})\not = \emptyset$$ for $i+j
\leq n$, $i, j = 1, 2, \cdots , n$ are contained in $C_0$ (see \cite{S} 
and \cite{H} for details on root compactness).  

Let ${\cal L}(E)$ be the space of all invertible continuous linear
operators on $E$: it may be noted that any $T\in {\cal L}(E)$ has 
continuous inverse (see 2.11 and 2.12 of \cite{Ru}).  
For $T \in {\cal L}(E)$ and $\mu \in {\cal P}(E)$,
$T(\mu )$ is defined by $T(\mu )(B) = \mu (T^{-1} (B))$ for all Borel
subsets $B$ of $E$.  A linear operator $T$ in ${\cal L}(E)$ is called 
{\it contraction} if $T^n (x) \ra 0$ for all $x \in E$.  

We equip ${\cal L}(E)$ with strong operator topology, that is, a net 
$(T_\ap )$ in ${\cal L}(E)$ converges (in the strong operator topology) to 
$T$ or $T_\ap \ra T$ 
if $T_\ap (x) \ra T(x)$ for all $x \in E$.   It may be noted that 
${\cal L}(E)$ equipped with strong topology and ${\cal P}(E)$ with weak 
topology, the map 
$\Psi \colon {\cal L}(E) \times {\cal P}(E) \ra {\cal P}(E)$ defined by
$\Psi (T, \lam ) =  T(\lam )$ is sequentially continuous.  

\bd\label{d1}
We say that a $\mu \in {\cal P}(E)$ is {\it $T$-semi-selfdecomposable} for 
a contraction $T \in {\cal L}(E)$ if there exists a sequence $(\nu _i)$ of 
probability measures on $E$ and an increasing sequence 
$N(n) \ra \infty $ and a sequence $(T_n)$
from ${\cal L}(E)$ such that $(\nu _{n,i})$ is a infinitesimal triangular
system with $\prod \nu _{n, i} x_n \ra \mu$ where 
$\nu _{n,i} = T_n (\nu_i)$ for $1\leq i \leq N(n)$, $(x_n)$ is a sequence 
in $E$ and $T_{n+1}T_n^{-1} \ra T$.
\ed

This notion is a natural extension of stable and semi-stable measures.  
The classical situation is studied in \cite{Lo}, for measures on
finite-dimensional vector spaces this notion was considered in 
\cite{MN} and \cite{MSW} and in \cite{Y} such measures on locally compact 
abelian groups 
was considered with emphasis on totally disconnected groups.

\bd\label{d2}
We call a measure $\mu \in {\cal P}(E)$ 
{\it $T$-decomposable} for a $T\in {\cal L}(E)$ if $T(\mu )$ is a factor 
of $\mu$, that is, $\mu = T(\mu ) \nu$ for a $\nu \in {\cal P}(E)$ (such a
$\nu$ is called corresponding co-factor) and in addition if 
$T^n(\mu ) \ra \delta _0$, $\mu$ is called {\it strongly $T$-decomposable}.  
\ed

The notion of $T$-decomposable measures on
infinite-dimensional cases, mainly the Banach space situation was studied
in \cite{Kr}, \cite{S0} and \cite{S1}.  

In this article we consider $T$-semi-selfdecomposable measures and prove 
that a probability measure $\mu $ is $T$-semi-selfdecomposable if and 
only if $\mu$ is (strongly) $T$-decomposable with an unique infinitely 
divisible co-factor.   

We next introduce subclasses of ${\cal P}(E)$ and $I(E)$ and study their 
basic properties that are similar to results in section 2 of \cite{MSW} 
(see also \cite{KS1} and \cite{KS2}) where the results are proved for 
measures on finite-dimensional vector spaces.  

A one-parameter subgroup $(T_t)$ of ${\cal L}(E)$ is a continuous 
homomorphism $t\mapsto T_t$ from $(0, \infty)$ into ${\cal L}(E)$.  
Let ${\cal L}_+$ be the set of all one-parameter subgroups $(T_t )$ such 
that $T_t (x) \ra 0$ as $t \ra 0$ for all $x \in E$.  

\bd\label{d3}
Let $0< b <1$, $(T_t) \in {\cal L}_+$ and $H \subset {\cal P}(E)$.  A 
measure
$\mu \in {\cal P}(E)$ is said to be in the class $\tilde K(H, b , (T_t))$ 
if
there exist a sequence of $(\mu _n)$ in $H$, a sequence $(a_n)$ of 
positive numbers and a subsequence $(k_n)$ of integers such that 
$${a_n \over a_{n+1}} \ra b ~~{\rm and }~~ T_{1\over a_n} 
(\mu _1 \mu _2 \cdots \mu _{k_n})x_n \ra \mu $$ for some sequence $(x_n)$ 
in $E$ and in addition 
if the system $\{ T_{1\over a_n} (\mu _i )\mid 1\leq i \leq k_n \}$
is infinitesimal, we say that $\mu \in K(H, b, (T_t))$.  
\ed

It follows from the definition that $K(H, b, (T_t)) \subset 
\tilde K(H, b, (T_t))$.  The class $\tilde K$ is studied in \cite{Bu} and 
the class $K$ is studied in \cite{MSW} for measures on finite-dimensional 
vector spaces.  

We now introduce other nested subclasses.  

\bd\label{d4}
For $0< b<1$ and $(T_t) \in {\cal L}_+$, we define 
$$L_0 (b, (T_t)) = K({\cal P}(E), b, (T_t)),$$
$$L_m (b, (T_t)) = K(L_{m-1}(b, (T_t)), b, (T_t)), m = 1,2, \cdots, $$ and 
$$L_\infty (b, (T_t)) = \cap _{m=1}^\infty L_m (b, (T_t)).$$ Similarly we 
define $\tilde L _m (b, (T_t) )$ using $\tilde K$ instead of $K$.  
\ed

\bd\label{d5}
For $(T_t) \in {\cal L}_+$ and $H\subset {\cal P}(E)$, we define 
$$K(H, [0, 1], (T_t)) = \cap _{b\in [0, 1]\setminus \{0,1 \}}
K(H, b ,(T_t))$$ and 
$$L_0([0,1], (T_t)) = K({\cal P}(E), [0,1], (T_t))$$ and 
we define inductively for $1\leq m < \infty$, 
$$L_m ([0,1], (T_t)) = K(L_{m-1}([0,1], (T_t)), [0,1], (T_t)),$$ 
$$L_\infty([0,1], (T_t)) = \cap L_m([0, 1], (T_t)).$$
Similarly define $\tilde L_m ([0,1], (T_t))$
using $\tilde K$ instead of $K$ for all $0\leq m \leq \infty$.
\ed

\bd\label{d6}
A $\mu \in {\cal P}(E)$ is called {\it operator selfdecomposable} 
with respect to $(T_t)\in {\cal L}^+$ if there
exist a sequence $(\mu _n)$ in ${\cal P}(E)$, a sequence $(x_n )$ in $E$ 
and $b_n >0$ such that the triangular system 
$\{ T_{1\over b_n}  (\mu _i) \mid 1\leq i \leq n \}$ is infinitesimal and 
$$T_{1\over b_n}(\prod _{j=1}^n\mu _j) x_n \ra \mu .$$  
\ed

The class of operator selfdecomposable measures was studied in \cite{MSW} 
for measures on finite-dimensional vector spaces and in \cite {J1}, 
\cite{KS1}, \cite{KS2}, \cite{Mi} and \cite{U1} this class was considered 
for measures on Banach spaces.  Here we prove the equivalence of operator 
selfdecomposable measures with other subclasses in Definition \ref{d5}.

We wish to remark that the operator selfdecomposable measures are also 
known as L\'evy's measures and the standard definition of such measures 
given by K. Urbanik \cite{U1} is different from the one given 
here.  In \cite{Mi} and \cite{U1}, the authors considered the group of 
invertible operators with norm topology and in particular, the 
semigroup $(T_t)$ is supposed to have bounded generator.  In 
\cite{KS1} and \cite{KS2}, the authors considered the group $(tI)$ of 
operators.  
Thus, when $E$ is a Banach space, Theorem \ref{th3} (to be proved in 
Section 4) and the main result in \cite{Mi} and \cite{U1} show that our 
notion is more general than the one given in \cite{Mi} and \cite{U1} 
provided the measure is not supported on a proper hyperplane.  
We also wish to mention the work of Riddhi Shah \cite{Rs} which further 
generalizes the notion of self-decomposability for non-commutative 
groups.  
\end{section}

\begin{section}{Preliminaries}

We first prove a technical lemma.  On complete separable metric
groups, it is known that point-wise contracting automorphisms are 
uniformly contracting on compact sets (see \cite{S2}).  Using an 
elementary 
argument we extend the result to 

\bl\label{lem1}
Let $T$ be a contracting operator in ${\cal L}(E)$.  Suppose $\Ga$ is
a relatively compact subset of ${\cal P}(E)$.  Then 
$T ^n (\lam ) \ra \delta _0$ uniformly for all $\lam \in \Ga$.
\el

\bo
Let $U$ be a neighborhood of $0$ in $E$ and $\epsilon >0$.  
Since $\Ga$
is relatively compact, there exists a compact set $K$ such that $\lam
(E\setminus K) <\epsi$ for all $\lam \in \Ga$.  By Lemma 1 of \cite{S2},
there exists a $N$ such that $T^n(K) \subset U$ for all $n \geq N$.  
Now for $n \geq N$, $T^n (\lam )(E\setminus U) \leq \lam (E\setminus K)
< \epsi$ for all $\lam \in \Ga$.  Thus, $T^n (\lam ) \ra \delta _0$
uniformly for all $\lam$ in $\Ga$.
\hfill{$\fbox{}$}
\eo

The proof of the following embeddability result is more or less identical 
to the proofs following W. B\"{o}ge as in \cite{S}. 

\bp\label{prop0}
Let $E$ be strongly root compact and $\mu \in {\cal P}(E)$.  Suppose
$(\mu _n)$ is a sequence in ${\cal P}(E)$ such that 
$\mu _n ^{k_n }\ra \mu$ as $k_n \ra \infty$.  Then $\mu$ is embeddable in
a c.c.s. $(\mu _t)$.
\ep

The {\it linear dual} $E^*$ of $E$ is the space of all continuous linear
functionals on $E$.  The following result may be seen as the generalized
central limit theorem on $E$ (see Chapter 2 of \cite{AG}).  

\bt\label{thm0}
Let $\mu \in {\cal P}(E)$.
Suppose $E$ is strongly root compact and the linear dual of $E$ separates 
points of $E$. Then the following are equivalent:  

\be
\item there exists an infinitesimal triangular system $\{ \mu _{n,i} \mid
1\leq i \leq k_n \}$ and a sequence $(x_n)$ in $E$ such that
$\prod _i \mu _{n, i} x_n \ra \mu$;

\item $\mu$ is infinitely divisible;

\item $\mu$ is embeddable in a c.c.s. $(\mu _t)$.

\ee
\et

\br
Any complete locally convex space with a countable base is a complete
separable
metric space which is strongly root compact (see \cite{S}) with separating
dual.  Also the spaces $L^p(\Omega, \Sigma, \lam)$ for any measure space 
$(\Omega, \Sigma, \lam)$ and for $0<p<1$ are complete separable metric 
spaces with separating dual.  It may also be proved using the arguments in 
\cite{BG} that the spaces $L^p(\Omega, \Sigma, \lam)$ are strongly root 
compact but are not locally convex. 
\er

\bo
We first show that (1) implies (2) by applying Hungarian semigroup theory 
of Rusza and Szekely (see \cite{RS}).    
Suppose $\mu$ is a limit of an infinitesimal triangular system.  
Since $E$ is a complete separable metric abelian group, ${\cal P}(E)$ is a 
stable Hungarian semigroup (see \cite{RS}).  We now claim that for any
$\lam \not = \delta _x$, there exists a map $f \colon F(\lam ) \ra [0,
\infty)$ such that $f(\lam _1 \lam _2) = f(\lam _1 ) +f(\lam _2)$ if $\lam
_1 , \lam _2 , \lam _1\lam _2 \in F(\lam )$ where $F(\lam )$ is the set of
factors of $\lam$.  Let $\lam $ be a non-degenerate probability measure on
$E$.  Then since the linear dual of $E$ separates
points of $E$, there exists a continuous linear functional $\phi$ such
that $\phi (\lam \check \lam ) \not = \delta _0$.  Then there exists a 
real number $r$
such that $$0< \int \exp (ir\phi (x)) d\lam \check \lam (x) <1$$ and hence 
$$ 0< \int \exp (ir\phi (x) ) d\lam _1 \check \lam _1 (x)$$ for any 
$\lam _1 \in F(\lam ) $.  Let $f \colon F(\lam ) \ra [0, \infty )$ be 
defined by $$f(\lam _1 ) = -\log (\int \exp (ir\phi (x)) d\lam _1 \check
\lam _1 (x) ) $$ for any $\lam _1 \in F(\lam )$.  This shows that ${\cal
P}(E)$ is a stable normable Hungarian semigroup.  By Theorem
26.9 of \cite{RS} and since any element of $E$ is infinitely divisible, 
$\mu$ is infinitely divisible.

Proposition \ref{prop0} shows that (2) implies (3).  We now prove
(3) implies (1).  Suppose $\mu$ is embeddable in a c.c.s. 
$(\mu _t)$.  Then for $n \geq 1$ and $1\leq i \leq n$, let 
$\mu _{n, i}=\mu _{{1\over n}}$.  Then $\prod _i \mu _{n,i} = \mu$ for
all $n \geq 1$.  We also have $\mu _{{1\over n}}\ra \delta _0$.  Thus,
$\mu$ is a limit of an infinitesimal triangular system.

\hfill{$\fbox{}$}
\eo

\br
\be
\item [(i)] Under the assumptions of Theorem \ref{thm0}, embedding c.c.s. 
is uniquely determined which may be seen as follows:  
Let $(\mu _t)$ and $(\lam _t)$ be two continuous convolution semigroups 
such that $\mu _1 = \lam _1$.  Since embedding c.c.s. is 
unique on $\mathbb R$, for any $f \in E^*$, $f (\lam _t) = f(\mu _t)$ for 
all $0 \leq t< \infty$ and hence $\lam _t = \mu _t$ for any $0<t<\infty$.  

\item [(ii)] Since any point measure $\delta _x$ is infinitely divisible, 
hence limit of an infinitesimal triangular array, condition (1) of 
Theorem \ref{thm0} is equivalent to 

\be
\item [(1')] there exists an infinitesimal triangular system 
$\{ \mu _{n,i} \mid 1\leq i \leq k_n \}$ such that
$\prod _i \mu _{n, i} \ra \mu$ (without shift term).
\ee
\ee
\er
\end{section}

\begin{section}{$T$-semi-selfdecomposable measures}

In this section we characterize $T$-semi-selfdecomposable measures on
vector spaces.  We first prove the following:

\bp\label{prp1}
Let $\mu$ be a probability measure on a complete separable metric 
vector space $E$.  Let $T$ be a contracting
operator in ${\cal L}(E)$ such that $\mu = T(\mu ) \nu$.  Suppose $\nu$ is
embeddable in a c.c.s $(\nu _t)$.  Then there exists a sequence $(\ga _i)$
in ${\cal P}(E)$ and an
increasing sequence $(N(n))$ of integers such that $N(n) \ra \infty$ and
the triangular system $\{ T^n (\ga _i) \mid 1\leq i \leq N(n), n\in
\mathbb N \}$ is infinitesimal and converges to $\mu$.
\ep

\bo
Since $\nu$ is embeddable in a c.c.s. $(\nu (t))$, let $\nu _n = \nu
({1\over n})  $ for all $n \geq 1$.  Let $N(n) = n^2$.  Then for each 
$n \geq 1$ and  $1\leq i \leq N(n)$ define 
$$\ga _i = T^{-k}\nu _{N(k) - N(k-1)}$$ if $N(k-1)<i\leq N(k) $.  
The assumuption $\mu = T(\mu )\nu$ implies $\mu = T^n(\mu ) \prod 
_{k=0}^{n-1} T^k(\nu )$ for any $n \geq 1$.  Then 

$$\begin{array}{rcl}
T^n (\ga _1\cdots \ga _{N(n)}) & = & T^n (\prod _{k=1}^n \prod
_{i=N(k-1)+1}^{N(k)} \ga _i ) \\ 
& = & T^n ( \prod _{k=1}^n T^{-k} \nu ) \\ 
& \ra & \mu
\end{array}$$
because $T^n (\mu ) \ra \delta _0$.  We now claim that the triangular
system $(T^n \ga _i)$ is infinitesimal.  Let $\Ga = \{ \nu (t) \mid
0\leq t\leq 1 \}$.  Then $\Ga$ is relatively compact.  Since $T$ is
contracting, by Lemma \ref{lem1} $T^n (\lam ) \ra \delta _0$ uniformly for
all $\lam \in \Ga$.  Let $U$ be a
neighborhood of $0$ in $E$ and $\epsilon >0$ be given.  Then there exists
a $N_0$ such that $T^n (\lam )(E\setminus U) < \epsi $ for all $n\geq
N_0$ and all $\lam \in \Ga$.  Now for $k \leq [{n\over 2}]$ and $n \geq
2N_0$, we have $n-k \geq N_0$ and hence 
$$\begin{array}{rcl}
T^n \ga _i (E\setminus U) & = & 
T^{n-k} \nu _{N(k)-N(k-1)} (E\setminus U)\\ 
& < & \epsilon 
\end{array}$$
for all $n \geq
2N_0$ and all $k \leq [{n\over 2}]$.  
Since $(\nu (t))$ is a c.c.s. we have $\nu _n
\ra \delta _0$.  Since $T$ is contracting, there exists a neighborhood
$V$ of $0$ such that $V \subset T^{-n} (U)$ for all $n \geq 1$.  For
$\epsilon >0$, there exists a $M_0$ such that $$\nu _n (E\setminus V)
<\epsilon $$ for all $n \geq M_0$.  Since $N(n) - N(n-1) \ra \infty$,
there exists a $N_1$ such that $N(n ) - N(n-1) > M_0$ for all $n \geq
N_1$.  Now, for $n \geq 2N_1$ and $k \geq [{n \over 2}]$ we have 
$N(k) - N(k-1) >M_0$.  Then we get
that $$\nu _{N(k)-N(k-1)} (E\setminus V) <\epsilon $$ for all 
$k\geq [{n\over 2}]$ and $n \geq 2N_1$.  This implies that 
$$\begin{array}{rcl}
T^n \ga _i(E\setminus U) & = &
T^{n-k}\nu _{N(k)-N(k-1)} (E\setminus U)\\
& \leq & \nu _{N(K)-N(k-1)} (E\setminus V) \\
& < & \epsilon
\end{array}$$
for all $k \geq [{n\over 2}]$ and $n \geq 2N_1$.  
Thus, $T^n \ga _i \ra \delta _0$ uniformly for all $\ga _i$, that is, for
every neighborhood $U$ of $0$ there exists $N>0$ such that 
$T^n (\ga _i)(E\setminus U) <\epsi$ for all $n \geq N$ and for all 
$i\geq 1$.  Thus, 
$(T^n (\nu _i ))$ is an infinitesimal triangular system converging to
$\mu$.  
\hfill{$\fbox{}$}
\eo

We now characterize $T$-semi-selfdecomposable measures on $E$
(see \cite{MN}, \cite{MSW} and section 3, Theorem 1 of
\cite{Y}).

\bt\label{thm1}
Let $E$ be a complete separable metric vector space which
is strongly root compact and the linear dual $E^*$ of $E$ separates points 
of $E$.  Let $\mu$ be a probability measure on $E$.  Let $T$ be a
contraction on $E$.  Then the following are equivalent:

\be

\item $\mu$ is $T$-semi-selfdecomposable;

\item $\mu$ is (strongly) $T$-decomposable with an infinitely divisible 
co-factor. 

In this case, the co-factors are uniquely determined.  
\ee

\et

\bo 
Let $$\lam _n = T_n (\nu _1 \cdots \nu _{N(n)})x_n ,$$ 
$$\rho _n = T_n (\nu _1 \cdots \nu _{N(n-1)})x_{n-1}$$ and 
$$\sigma _n = T_n (\nu _{N(n-1)+1} \cdots \nu _{N(n)})x_nx_{n-1}^{-1}$$
for all $n \geq 1$.  Since $\lam _n \ra \mu$, $\rho _n \ra T(\mu )$ 
and hence $\lam _n = \rho _n \sigma _n$ implies that $(\sigma _n )$ is
relatively compact.  Thus, for any limit
point $\nu$ of $(\sigma _n) $, we have $\mu = T(\mu ) \nu$.   Since 
$\nu$ is the limit of a infinitesimal triangular system, by Theorem
\ref{thm0}, $\nu$ is also infinitely divisible.  Since $T$ is contracting,
any $T$-decomposable measure is obviously strongly $T$-decomposable, 
hence $\mu$ is strongly $T$-decomposable.  Now for the uniqueness, 
let $\nu '\in {\cal P}(E)$ be such that $\mu = T(\mu ) \nu '$.  
By Theorem \ref{thm0} $\mu$ is infinitely divisible, hence $T(\mu )$ is 
also infinitely divisible.  Now for any $f \in E^*$,  
the characteristic function of $f(T(\mu ) )$ does not vanish and hence 
$f(\nu ) = f(\nu ')$.  Thus, $f(\nu ) = f(\nu ')$ for 
any $f \in E^*$.  This implies that $\nu = \nu '$.  

Conversely, suppose $\mu$ is $T$-decomposable with an
infinitely divisible co-factor, say $\nu$.  Then by Theorem \ref{thm0}
$\nu$ is embeddable in a c.c.s. $(\nu (t))$.  Now the result
follows from Proposition \ref{prp1}.   
\hfill{$\fbox{}$}
\eo
\end{section}

\begin{section}{Nested subclasses of ${\cal P}(E)$ and $I(E)$}

In this section we assume that $E$ is a strongly root
compact complete separable metric vector space and the linear dual of $E$
separates points of $E$ and we prove basic properties of the 
various subclasses defined in the introduction and their connection with 
operator selfdecomposable measures.  We start with the following notions. 

\bd
A subset $H\subset {\cal P}(E)$ is said to be $(T_t)$-completely closed if 
$H$ is a closed subsemigroup of $ {\cal P}(E)$ and $T_a\mu 
*\delta _x \in H$ whenever $a >0$, $x \in E$ and $\mu \in H$.  
When $H\subset I(E)$, we say that $H$ is $(T_t)$-completely closed in the 
strong sense if $H$ is $(T_t)$-completely closed and $\mu $ in $H$ is 
embeddable in an one-parameter subgroup in $H$.
\ed  

It is easy to see that $I(E)$ is $(T_t)$-completely closed in the strong 
sense.  We first consider the classes $K$ and $\tilde K$.  

\bt\label{th1}

Let $H \subset {\cal P}(E)$ and $(T_t)\in {\cal L}_+$.  Then we have the
following:

\be

\item $K(H, b, (T_t)) \subset I(E)$; 

\item Suppose $H$ is $(T_t)$-completely closed.  If 
$\mu \in K(H, b, (T_t))$, then $\mu$ is (strongly) $T_b$-decomposable 
with a unique co-factor $\mu _b \in H\cap I(E)$; 

\item Suppose $H$ is $(T_t)$-completely closed in the strong sense, then 
the converse of (2) is also true and $K(H, b, (T_t))$ is also 
$(T_t)$-completely closed in the strong sense.
\ee
\et

\bo
It follows from Theorem \ref{thm0} that $K(H, b, (T_t)) \subset I(E)$.  
We next observe that if ${a_n \over a_{n+1}} \ra b$, then
$T_{1\over a_{n+1}} (T_{1\over a_n })^{-1} \ra T_b$.  Then statement (2) 
and the first part of statement (3) follow from 
Theorem \ref{thm1} and the second part of statement (3) can be easily 
verified.  
\hfill{$\fbox{}$}
\eo

In a similar way we may prove 

\bt\label{th2}

Let $H \subset {\cal P}(E)$ and $(T_t) \in {\cal L}_+$.  Suppose $H$ is
$(T_t)$-completely closed.  Then we have the following.

\be

\item $\tilde K(H, b, (T_t))$ is also $(T_t)$-completely closed;

\item $\mu \in \tilde  K(H, b, (T_t))$ if and only if $\mu$ is (strongly) 
$T_b$-decomposable with a (not necessarily unique) co-factor 
$\mu _b \in H$.

\ee
\et

It follows from Theorem \ref{thm0} that $(L_m (b, (T_t)))$ are decreasing 
subclasses of $I(E)$.  We next consider the nested subclasses $(L_m)$ and 
$(\tilde L_m)$.  

\bp\label{p1}
For $0\leq m \leq \infty$, we have the following:
\be
\item $\mu \in L_m (b, (T_t))$ if and only if $\mu$ is (strongly)
$T_b$-decomposable with a unique co-factor $\rho \in L_{m-1}(b, (T_t))$ 
where $L_{-1}$ and $L_{\infty -1}$ are understood as $I(E)$ and 
$L_\infty$.  Also, $L_m (b, (T_t))$ is $(T_t)$-completely closed in the 
strong sense;

\item $\mu \in \tilde L_m (b, (T_t))$ if and only if $\mu$ is (strongly)
$T_b$-decomposable with a (not necessarily unique) co-factor $\rho _m \in 
\tilde L_{m-1}(b, (T_t))$ where $\tilde L_{-1}$ and $\tilde L_{\infty -1}$ 
are understood as 
${\cal P} (E)$ and $\tilde L_\infty $.  Also, $\tilde L_m (b, (T_t))$
is $(T_t)$-completely closed.

\ee
\ep

\bo
It may be easily seen that the first statement follows from Theorem 
\ref{th1} and since $I(E)$ is $(T_t)$-completely closed in the strong 
sense.   

Now for the second part, since ${\cal P}(E)$ is completely closed, the 
second part for $m =0$, 
follows from Theorem \ref{th2}.  Now the second part follows by induction 
on $m$ for $1\leq m < \infty$.  For $m =\infty$, if $\mu \in \tilde 
L_\infty (b, (T_t))$, then for each $m < \infty$, there exists a 
$\rho _m\in \tilde L_m(b, (T_t))$ such that $\mu = T_b(\mu ) \rho _m$.  
This shows that $(\rho _m)$ is relatively compact.  There exists a subsequence
$(k_m)$ such that $\rho _{k_m} \ra \rho$.  This implies that $\mu =
T_b(\mu ) \rho$.  Also, since $\tilde L_m (b, (T_t)) \subset \tilde L_n(b, 
(T_t))$ if $n \leq m$, we get that $\rho \in \tilde L_m (b, (T_t))$ for 
all $m \geq 1$.  Thus, proving the second part.  
\hfill{$\fbox{}$}
\eo

Note that in view of Proposition \ref{p1}, for $0<b<1$ and 
$(T_t) \in {\cal L}_+$, $\mu \in L_0 (b, (T_t))$ if and only if $\mu$ is 
$T_b$-semi-selfdecomposable in the sense of Definition 1.

We next state another technical lemma which may be proved using 
Appendix 2 in \cite{S2} arguing as in Lemma \ref{lem1}.  

\bl\label{lm1}
Let $(T_t) \in {\cal L}_+$ and $\Ga$ be a relatively compact subset of 
${\cal P}(E)$.  Then as $t \ra 0$, $T_t (\lam ) \ra 0$
uniformly for all $\lam \in \Ga$, that is for given $\epsi >0$ and
a neighborhood $U$ of $0$ there exists a $s >0$ such that 
$T_t(\lam )(E\setminus U) < \epsi$ for all $t<s$ and all $\lam \in \Ga$.
\el

We now obtain the equivalence between operator selfdecomposable measures
and other subclasses defined in the introduction (Definition 6).   

\bt\label{th3}
Let $(T_t) \in {\cal L}_+$ and $\mu \in {\cal P}(E)$.  Then the following 
are equivalent:
\be

\item $\mu$ is operator selfdecomposable with respect to $(T_t)$;

\item $\mu \in L_0([0,1], (T_t))$;

\item $\mu \in \tilde L_0 ([0, 1], (T_t))$.

\ee
\et

\bo 
It is easy to see that (2) implies (3).  We now prove (1) implies (2).
Suppose $\mu$ is operator selfdecomposable. If $\mu = \delta _x$ for
some $x \in E$, then (1) implies (2).  So, we assume $\mu \not =
\delta _x$ for any $x \in E$.  Let $(b_n)$ be a sequence of positive reals 
and $(\mu _n)$ be a sequence of probability measures on $E$ such that 
$\{ T_{1\over b_n}(\mu _i ) \mid n \geq 1, ~~1\leq i \leq n \}$ is a 
infinitesimal triangular system and $T_{1\over b_n} (\mu _1 \mu _2 \cdots 
\mu _n )x _n \ra \mu $ for some sequence $(x_n)$ in $E$.  

We first claim that $b_n \ra \infty$ and 
${b_n \over b_{n+1}}\ra 1$.  Suppose $b_n \not \ra \infty$.  
Then there exists a subseqeunce $(b_{m(n)})$ of $(b_n)$ such 
that $b_{m(n)} \ra b <\infty$.  Since $\{ T_{1\over b_n}(\mu _i ) 
\mid n \geq 1, ~~1\leq i \leq n \}$ is infinitesimal, $T_{1\over b_n} (\mu 
_i) \ra \delta _0$ for all $i\geq 1$.  This implies that $T_{1\over b} 
(\mu _i) = \delta _0$ for all $i \geq 1$.  Since $T_{1\over b}$ is 
invertible, $\mu _i = \delta _0$ for all $i \geq 1$.  This shows that 
$\mu$ itself is a dirac measure which is a contradiction.  Thus, $b_n \ra 
\infty$.  For any $\nu \in {\cal P}(E)$, let $\check \nu \in {\cal 
P}(E)$ be defined by $\check \nu (B) = \nu (-B)$ for any Borel subset 
$B$ of $E$.  Now let $\mu ' = \mu \check \mu$ and $\mu _n ' = \mu _n 
\check \mu _n$ for any $n \geq 1$.  Then $T_{1\over b_n} (\mu _1' \mu _2' 
\cdots \mu _n ' ) \ra \mu '$.  Since $T_{1\over b_n} (\mu _n) \ra \delta 
_0$, we have $T_{1\over b_{n+1}} (\mu _1' \mu _2' \cdots \mu _n' )\ra \mu 
'$.  For any $n \geq 1$, let $\rho _n = T_{1\over b_n} (\mu _1' \mu _2' 
\cdots \mu _n')$.  Then $T_{b_n\over b_{n+1}} (\rho _n) \ra \mu '$.  
Suppose $({b_n \over b_{n+1}})$ is not bounded, then there exists a 
subsequence $k(n)$ such that ${b_{k(n)}\over b_{k(n)+1}} \ra \infty$.  Let 
$a_n = {b_{k(n)}\over b_{k(n)+1}}$.  Then by Lemma \ref{lm1}, $T_{1\over 
a_n} (T_{b_m\over b_{m+1}}(\rho _m) )\ra \delta _0$ uniformly for all 
$m \geq 1$.  Let $U$ be a neighbourhood of $0$ and $\epsi >0$.  Then there 
exists a $N$ such that $T_{1\over a_n} (T_{b_m\over b_{m+1}}(\rho _m) ) 
(U) >1-\epsi$ for all $n \geq N$ and all $m \geq 1$.  Now for $n \geq N$, 
$T_{b_{k(n)+1}\over b_{k(n)}} T_{b_{k(n)} \over b_{k(n)+1}} 
(\rho _{k(n)}) (U) >1-\epsi$.  
This implies that $\rho _{k(n)}(U) >1-\epsi$ for all $n 
\geq N$, that is, $\rho _{k(n)} \ra \delta _0$ but $\rho _n \ra \mu '$.  
Thus, $\mu ' = \delta _0$ and hence $\mu$ is a dirac measure which is a 
contradiction.  Thus, $({b_n \over b_{n+1}})$ is bounded.  Suppose $b$ is 
a limit point of $({b_n \over b_{n+1}})$.  
Then since $T_{b_n\over b_{n+1}} (\rho _n ) \ra \mu '$ and 
$\rho _n\ra \mu '$, we have 
$T_b(\mu ' ) = \mu '$.  If $b\not = 1$, then $\mu ' =T_b^n (\mu ') \ra 
\delta _0$ or $\mu' = T_b^{-n} (\mu ') \ra \delta _0$ and therefore $\mu$ 
is a dirac measure.  This is a contradiction and hence $b=1$.  Thus, 
${b_n \over b_{n+1}} \ra 1$.  

Now for any $0< b <1$, there exist sequences
$(m_k)$ and $(n _k)$ such that $m _k < n _k$ such that ${b_{m _k} \over
b_{n_k}}\ra b$.  Now for each $k \geq 1$, we have 
$$T_{1\over b_{n_k}} T_{b_{m_k}} T_{1\over b_{m_k}} ( \mu _1 \cdots \mu 
_{m_k}) x_{m_k}
T_{1\over b_{n_k}} (\mu _{m _k+1} \cdots \mu _{n_k})x_{n_k}x_{m_k}^{-1}\ra
\mu$$ and hence 
$$\mu = T_b(\mu ) \rho$$ where $\rho$ is a limit point of 
$T_{1\over b_{n_k}} (\mu _{m_k+1} \cdots \mu _{n_k})x_{n_k}x_{m_k}^{-1}$.  
This shows that $\rho$ is infinitely divisible (by Theorem \ref{thm0} ).

We now prove (3) implies (1).  Suppose $\mu \in \tilde 
L([0, 1],(T_t))$.  Now for each $n \geq 1$, 
let $a_n = {n\over n+1}$.  
Then there exists a $\mu _n \in {\cal P}(E)$ such that 
$\mu = T_{a_n} (\mu ) \mu _n$.  Let $b_n = {1\over n}$ and 
$\nu _n = T_{n+1} (\mu _n)$.  Then 
$T_{n+1}(\mu ) =  T_{n+1}T_{a_n} (\mu ) \nu _n = 
T_{n}(\mu ) \nu _n$ for all $n \geq 1$.  This implies that 

$$\begin{array}{rcl}
T_{b_n}  (\nu _1 \nu _2 \cdots \nu _n  ) T_{b_n}(\mu ) 
& = &  T_{b_n} (\mu \nu _1 \cdots \nu _n) \\
 & = & T_{b_n}( T_2(\mu ) \nu _2 \cdots \nu _n) \\
& = & T_{b_n} (T_{n+1} (\mu ))\\
& = & T_{1\over a_n} (\mu )
\end{array}$$
for all $n \geq 1$.  Since ${n+1 \over n} \ra 1$ and $b_n\ra 0$, we get
that $T_{b_n}(\nu _1 \nu _2 \cdots \nu _n ) \ra \mu$.  
We now prove that $\{ T_{b_n}(\nu _j) \mid 1\leq j \leq n,~ n \geq 1 \}$ 
is an infinitesimal triangular system.  Now, for $1\leq j \leq n$, 
$T_{b_n} (\nu _j) = T_{j+1\over n}(\mu _j )$.  For $1\leq j \leq 
\sqrt n$, ${j+1\over n} \ra 0$.  Also, Since $\mu = T_{a_n} (\mu )
\mu _n$ and $a _n \ra 1$, we get that $\mu _n\ra \delta _0$.  Thus, 
by Lemma \ref{lm1}, for $\epsi >0$ and for a neighborhood $U$ of $0$,
there exists a $M_1$ such that $T_{b_n}(\nu _j ) (E\setminus U) < \epsi $
for all $n \geq M_1$ and all $j \leq {\sqrt n}$.  Now for all $1 \leq j
\leq n$,  ${j+1\over n} \leq 2$.  Since $\{ T_t (x) \mid 0 < t \leq 2 \}$
is relatively compact, by uniform boundedness principle, $\{ T_t \mid 0< t
\leq 2 \}$ is equicontinuous.  Let $U$ be any neighborhood of $0$ and
$\epsi >0$.  Let $V$ be a neighborhood of $0$ such that $T_t(V)
\subset U$ for all $0<t\leq 2$.  
Since $\mu _n \ra \delta _0$, there exists a $M_2$ such that
$\mu _n (E\setminus V) <  \epsi$ for all $n \geq M_2$.  This implies that 
$T_t(\mu _n (E\setminus U) ) \leq \mu _n (E \setminus V) < \epsi$ for all
$n \geq M_2$ and for all $0< t \leq 2$.  Thus, for $ n \geq M_2 ^2$
and ${\sqrt n}\leq j \leq n$, $T_{b_n} (\nu _j )(E\setminus U) \leq
T_{j+1\over n} (\mu _j)(E\setminus U) < \epsi$.  This shows that 
$T_{b_n}(\nu _j ) \ra \delta _0$ uniformly for all $1\leq j \leq n$, that
is $\{ T_{b_n} (\nu _j) \}$ is an infinitesimal triangular system.   
\hfill{$\fbox{}$}
\eo

\br
In the finite-dimensional situation, Theorem \ref{th3} was stated as 
Proposition 2.5 in \cite{MSW} and the proof was refered to Theorem 2.1 and 
Corollary 2.4 of \cite{Sa} where the result was proved for the 
one-parameter group $(tI)$ of homothetical operators and to 
Theorem 3.3.5 of \cite{JM} where the result was proved for full 
measures.  In view of this the proof of Proposition 2.5 in \cite{MSW} is 
not correct.  Thus, Theorem \ref{th3} corrects the proof of 
Proposition 2.5 of \cite{MSW} and extends it to any complete separable 
vector space.  
\er

The following results can be proved as in Theorem 2.2 of \cite{MSW} 
and we omit the details.

\bt\label{th4}
Let $E$, $\mu$ and $(T_t)$ be as in Theorem \ref{th3}.  Then for 
$0\leq m \leq \infty$, 
\be
\item $L_m ([0, 1], (T_t))= \cap _{b\in (0, 1)}L_m (b, (T_t))$;

\item $\tilde L_m ([0, 1], (T_t)) = \cap _{b\in (0, 1)} \tilde L_m (b, 
(T_t))$;

\item $L_m ([0, 1], (T_t)) = \tilde L_m ([0, 1], (T_t))$.
\ee

\et

\end{section}

\begin{acknowledgement}
I would like to thank the referee for his suggestions in improving the
text. 
\end{acknowledgement}

\vskip 0.25in

\noindent {C. Robinson Edward Raja, \\
Indian Statistical Institute,\\
Statistics and Mathematics Unit,\\
8th Mile Mysore Road,\\
R. V. College Post,\\
Bangalore - 560 059.\\
India.}

\noindent {creraja@isibang.ac.in}

\end{document}